**Some elementary remarks on the powers of a partial theta function**

**and corresponding q-analogs of the binomial coefficients.**

Johann Cigler

**Abstract**

We obtain formulas for the coefficients of positive and negative powers of a partial theta function.

## 1. Introduction

The partial theta function

(1)
$$f(x)=\sum_{n\ge0} q^{\binom{n+1}{2}} x^n$$

can be interpreted as a $q-$analog of the geometric series $\frac{1}{1-x}=\sum_{n\ge0} x^n$.

Since $\frac{1}{(1-x)^{k+1}}=\sum_{n\ge0}\binom{n+k}{k}x^n$ the coefficients of $f(x)^{k+1}$ are $q-$analogs of the binomial coefficients $\binom{n+k}{k}$. If we write

(2)
$$f(x)^{k+1}=\sum_{n\ge0}\left\langle {n+k \atop k} \right\rangle_q x^n$$

we get a $q-$Pascal triangle $\left(\left\langle {n \atop k}\right\rangle_q\right)_{n,k\ge0}$. The first terms are

(3)
$$\left(\left\langle {n \atop k}\right\rangle_q\right)_{n,k=0}^5=\begin{pmatrix} 1 & 0 & 0 & 0 & 0 & 0\\ q & 1 & 0 & 0 & 0 & 0\\ q^3 & 2q & 1 & 0 & 0 & 0\\ q^6 & q^2(1+2q) & 3q & 1 & 0 & 0\\ q^{10} & 2q^4(1+q^2) & 3q^2(1+q) & 4q & 1 & 0\\ q^{15} & q^6(1+2q+2q^4) & q^3(1+6q+3q^3) & 2q^2(3+2q) & 5q & 1 \end{pmatrix}$$

We first derive some properties and a combinatorial interpretation of these $q-$binomial coefficients. Then we determine the coefficients of the series $\frac{1}{f(x)^k}$.

## 2. Some properties of these q-binomial coefficients

**2.1.** The identity $f(x)^{k+1}=f(x)f(x)^k$ gives the recursion

(4) $$\left\langle {n \atop k} \right\rangle_q = \sum_{j=0}^{n-k} q^{\binom{j+1}{2}} \left\langle {n-j-1 \atop k-1} \right\rangle_q$$

with $\left\langle {n \atop 0} \right\rangle_q = q^{\binom{n+1}{2}}$ and $\left\langle {0 \atop k} \right\rangle_q = [k=0].$

For example,

$$\left\langle {5 \atop 2} \right\rangle_q = \left\langle {4 \atop 1} \right\rangle_q + q\left\langle {3 \atop 1} \right\rangle_q + q^3\left\langle {2 \atop 1} \right\rangle_q + q^6\left\langle {1 \atop 1} \right\rangle_q = 2q^4\left(1+q^2\right)+q\left(q^2\left(1+2q\right)\right)+q^3\left(2q\right)+q^6 = q^3+6q^4+3q^6.$$

**2.2.** Since $q^{\binom{n+1}{2}} = q^{\binom{n}{2}}q^n$ we see that $f(x)=1+qxf(qx)$ and $f(x)^{k+1}=\left(1+qxf(qx)\right)^{k+1}=\sum_{j=0}^{k+1}\binom{k+1}{j}q^jx^jf(qx)^j$ which implies

(5) $$\left\langle {n \atop k} \right\rangle_q = q^{n-k}\sum_{j=0}^{k+1}\binom{k+1}{j}\left\langle {n-k-1 \atop j-1} \right\rangle_q$$

for $0 \le k \le n$, where on the right side $\left\langle {-1 \atop j} \right\rangle_q = 0$ for $j \in \mathbb{N}$, $\left\langle {-1 \atop -1} \right\rangle_q = 1$ and $\left\langle {n \atop -1} \right\rangle_q = 0$ for $n \in \mathbb{N}$.

For example,

$$\left\langle {5 \atop 2} \right\rangle_q = q^3\left(3\left\langle {2 \atop 0} \right\rangle_q + 3\left\langle {2 \atop 1} \right\rangle_q + \left\langle {2 \atop 2} \right\rangle_q\right) = q^3\left(3q^3+3(2q)+1\right) = q^3+6q^4+3q^6.$$

(5) also implies that $\left\langle {n \atop k} \right\rangle_q \in \mathbb{N}[q]$ with highest term $(k+1)q^{\binom{n-k+1}{2}}$ for $0 \le k < n$.

### 2.3. A combinatorial interpretation

**Theorem 1**

*The $q-$binomial coefficients $\left\langle {n \atop k} \right\rangle_q$ have the following combinatorial interpretation:*

*If we define the weight of an integer partition $\lambda=\left(\lambda_1,\lambda_2,\cdots,\lambda_\ell\right)$ of $n=\lambda_1+\lambda_2+\cdots+\lambda_\ell$ with $\lambda_1 \ge \lambda_2 \ge \cdots \ge \lambda_\ell$ by $w(\lambda)=coef\left(\lambda\right)q^{ex(\lambda)}$ with $coef(\lambda)=\prod_{i=1}^{\ell-1}\binom{\lambda_i}{\lambda_{i+1}}$ and $ex(\lambda)=\sum_{i=1}^{\ell}(i-1)\lambda_i$*

*then*

(6) $$\left\langle {n \atop k} \right\rangle_q = \sum_{\lambda \in P_{n+1,k+1}} w(\lambda)$$

*where* $P_{n.k}$ *denotes the set of all partitions of* $n$ *with first term* $\lambda_1 = k$.

For example let us consider $\left\langle {5 \atop 2} \right\rangle_q$ :

$P_{6,3} = \{(3,3),(3,2,1),(3,1,1,1)\}$ implies

$$\left\langle {5 \atop 2} \right\rangle_q = w(3,3) + w(3,2,1) + w[3,1,1,1) = \binom{3}{3} q^3 + \binom{3}{2}\binom{2}{1} q^{2+2} + \binom{3}{1} q^{1+2+3} = q^3 + 6q^4 + 3q^6.$$

**Remark**

I want to thank Michael Schlosser for the hint to [2] and conjecturing formula (6). In [2] Peter Luschny stated without proof that $\sum_{\lambda \in P_{n,k}} coef(\lambda) = \binom{n-1}{k-1}$. Since I did not see how to prove this assertion I posted question [1]. As a comment Marko Riedel [3] gave a link to his answer of question MSE 4940765. His ideas led to the following

**Proof of Theorem 1**

To each $\lambda = (\lambda_1, \lambda_2, \cdots, \lambda_\ell) = (k, \lambda_2, \cdots, \lambda_\ell) \in P_{n,k}$ we associate $i(\lambda) = (i_1, i_2, \cdots, i_\ell, 0, \cdots, 0) \in \mathbb{N}^n$ with $i_j = \lambda_j - \lambda_{j+1}$ for $1 \le j < \ell$, $i_\ell = \lambda_\ell$ and $i_j = 0$ for $j > \ell$.

Then $\sum_{j=1}^{n} i_j = k$ and since $\lambda_j = i_j + i_{j+1} + \cdots + i_\ell$ we have

$$n = \lambda_1 + \lambda_2 + \cdots + \lambda_\ell = i_\ell + (i_\ell + i_{\ell-1}) + (i_\ell + i_{\ell-1} + i_{\ell-2}) + \cdots = \ell i_\ell + (\ell - 1) i_{\ell-1} + \cdots + i_1 = \sum_{j=1}^{n} j i_j.$$

On the other hand to each $i = (i_j)_{j=1}^{n}$ with $\sum_{j=1}^{n} i_j = k$ and $\sum_{j=1}^{n} j i_j = n$ there exists a uniquely defined partition $\lambda \in P_{n,k}$.

We then get

(7) $$w(\lambda) = \frac{k!}{\prod_{j=1}^{n} i_j!} q^{\sum_{j=1}^{n} \binom{j}{2} i_j}$$

since

$$\sum_{j=1}^{\ell}(j-1)\lambda_j = \left(i_2+i_3+i_4+\cdots\right)+2\left(i_3+i_4+i_5+\cdots\right)+\cdots+(\ell-1)i_\ell = \binom{1}{2}i_1+\binom{2}{2}i_2+\cdots+\binom{\ell}{2}i_\ell$$
$$=\sum_{j=1}^{n}\binom{j}{2}i_j$$

and

$$w(\lambda)=\prod_{i=1}^{\ell(\lambda)}\binom{\lambda_i}{\lambda_{i+1}}q^{\sum_{i=1}^{\ell}(i-1)\lambda_i} = \prod_{i=1}^{\ell(\lambda)-1}\frac{\lambda_i!}{\lambda_{i+1}!(\lambda_i-\lambda_{i+1})!}q^{\sum_{i=1}^{\ell}(i-1)\lambda_i} = \frac{k!}{\lambda_\ell!(\lambda_{\ell-1}-\lambda_\ell)!(\lambda_{\ell-2}-\lambda_{\ell-1})!\cdots(\lambda_1-\lambda_2)!}q^{\sum_{i=1}^{\ell}(i-1)\lambda_i}$$
$$=\frac{k!}{\prod_{j=1}^{n}i_j!}q^{\sum_{j=1}^{n}\binom{j}{2}i_j}.$$

The multinomial theorem

$$\left(x_1+x_2+\cdots+x_\ell\right)^{k+1}=\sum_{i_1+i_2+\cdots+i_\ell=k+1}\frac{(k+1)!}{i_1!i_2!\cdots i_\ell!}x_1^{i_1}x_2^{i_2}\cdots x_\ell^{i_\ell}$$

gives

$$\left(q^{\binom{1}{2}}z+q^{\binom{2}{2}}z^2+\cdots+q^{\binom{n+1}{2}}z^{n+1}\right)^{k+1} = \sum_{i_1+\cdots+i_{n+1}=k+1}\frac{(k+1)!}{i_1!\cdots i_{n+1}!}q^{\binom{1}{2}i_1+\binom{2}{2}i_2+\cdots+\binom{n+1}{2}i_{n+1}}z^{i_1+2i_2+\cdots+(n+1)i_{n+1}}$$

Therefore,

(8) $$[z^{n+1}]\left(q^{\binom{1}{2}}z+q^{\binom{2}{2}}z^2+\cdots+q^{\binom{n+1}{2}}z^{n+1}\right)^{k+1} = \sum_{\lambda\in P_{n+1,k+1}}w(\lambda).$$

On the other hand we have

$$[z^{n+1}]\left(q^{\binom{1}{2}}z+q^{\binom{2}{2}}z^2+\cdots+q^{\binom{n+1}{2}}z^{n+1}\right)^{k+1} = [z^{n-k}]\left(1+qz+\cdots+q^{\binom{n+1}{2}}z^n\right)^{k+1}$$
$$=\left[z^{n-k}\right]f(z)^{k+1}=\left\langle {n \atop k}\right\rangle_q .$$

For example, for $(n+1,k+1)=(6,3)$ the $6-$tuples $\left(i_1,i_2,\cdots,i_6\right)$ satisfying $\sum_j i_j=3$ and $\sum_{j=1}^{6}ji_j=6$ are $\left\{\left(0,3,0,0,0,0\right),\left(1,1,1,0,0,0\right),\left(2,0,0,1,0,0\right)\right\}$.

Therefore

$$[z^6]\left(q^{\binom{1}{2}}z+q^{\binom{2}{2}}z^2+\cdots+q^{\binom{6}{2}}z^6\right)^3=\frac{3!}{3!}\left(q^{\binom{2}{2}}\right)^3+\frac{3!}{1!^3}\left(q^{\binom{1}{2}}\right)\left(q^{\binom{2}{2}}\right)\left(q^{\binom{3}{2}}\right)+\frac{3!}{2!1!}\left(q^{\binom{1}{2}}\right)^2\left(q^{\binom{4}{2}}\right)$$

$$=\left(q^3+3!q^{0+1+3}+\frac{3!}{2!}q^{0+6}\right)=\left(q^3+6q^4+3q^6\right)=\left\langle {5 \atop 2} \right\rangle_q .$$

## 3. Negative powers of $f(x)$.

**Theorem 2**

*Let* $\dfrac{1}{f(x)^k}=\sum_{n\geq 0}u_k(n,q)x^n$*. Then*

$$u_k(n,q)=-q^n\sum_{j=0}^{n-1}(-1)^j\left\langle {n-1 \atop j} \right\rangle_q\binom{k+j}{k-1}. \tag{9}$$

**Proof**

$$\frac{1}{f(x)^k}=\frac{1}{\left(1+qxf(qx)\right)^k}=\sum_{j\geq 0}\binom{j+k-1}{k-1}(-1)^j q^j x^j\sum_{\ell\geq 0}\left\langle {\ell+j-1 \atop j-1} \right\rangle_q q^\ell x^\ell$$

$$=\sum_{n\geq 0}x^n q^n\sum_{j+\ell=n}(-1)^j\binom{j+k-1}{k-1}\left\langle {\ell+j-1 \atop j-1} \right\rangle_q=\sum_{n\geq 0}x^n q^n\sum_{j=0}^{n-1}(-1)^{j-1}\left\langle {n-1 \atop j} \right\rangle_q\binom{k+j}{k-1}.$$

For $k=1$ writing $u_1(n,q)=u(n,q)$ we get

$$u(n,q)=-q^n\sum_{j=0}^{n-1}(-1)^j\left\langle {n-1 \atop j} \right\rangle_q . \tag{10}$$

For example, $u(1,q)=-q(\left\langle {0 \atop 0} \right\rangle_q)=-q,\ u(2,q)=-q^2\left(\left\langle {1 \atop 0} \right\rangle_q-\left\langle {1 \atop 1} \right\rangle_q\right)=-q^2\left(q-1\right)=q^2-q^3,$

$$u(3,q)=-q^3\left(\left\langle {2 \atop 0} \right\rangle_q-\left\langle {2 \atop 1} \right\rangle_q+\left\langle {2 \atop 2} \right\rangle_q\right)=-q^3\left(q^3-2q+1\right)=-q^3+2q^4-q^6.$$